\documentclass[A3paper,fleqn,showlink]{cas-sc}
\newdefinition{theorem}{Theorem}[section]
\newdefinition{lemma}[theorem]{Lemma}
\newdefinition{corollary}[theorem]{Corollary}

\newdefinition{remark}{Remark}[section]
\numberwithin{equation}{section}

\usepackage{lipsum}
\usepackage[bookmarks,A3paper,colorlinks=true]{}
\usepackage[]{natbib}

\usepackage{amsmath}
\usepackage{amssymb}
\usepackage{amsthm}
\usepackage{latexsym}
\usepackage{graphicx}

\xspaceaddexceptions{]}
\def\tsc#1{\csdef{#1}{\textsc{\lowercase{#1}}\xspace}}
\tsc{WGM}
\tsc{QE}
\tsc{EP}
\tsc{PMS}
\tsc{BEC}
\tsc{DE}
\pdfoutput=1
\begin{document}
\let\WriteBookmarks\relax
\def\floatpagepagefraction{1}
\def\textpagefraction{.001}
\shorttitle{Scalar BSDEs with integrable terminal values and monotonic generators}
\shortauthors{Hun O et al.}
%\begin{frontmatter}

\title [mode = title]{Uniqueness, Comparisons and Stability for Scalar BSDEs with $L\exp(\mu\sqrt{2\log(1+L)})$ -integrable terminal values and monotonic generators}                      

\author[1]{Hun O}
\author[1]{Mun- Chol Kim}
\author[1]{Chol- Kyu Pak}\cormark[1]
\address[1]{Faculty of Mathematics, Kim Il Sung University, Pyongyang, Democratic People's Republic of Korea}
\ead{*pck2016217@gmail.com}

\begin{abstract}
This paper considers a class of scalar backward stochastic differential equations (BSDEs) with $L\exp(\mu\sqrt{2\log(1+L)})$-integrable terminal values. We associate these BSDEs with BSDEs with integrable parameters through Girsanov change. Using this technique, we prove uniqueness, comparisons and stability for them under an extended monotonicity condition (more precisely one sided Osgood condition).\end{abstract}

\maketitle

\section{Introduction}
\par Let $(\Omega,\mathcal{F},\mathbb{P})$ be a probability space, $T>0$ a finite time and $W$ a standard $d$-dimensional Brownian motion. Let $\mathbb{F}:=\{\mathcal{F}_t\}_{0\le t\le T}$ be a completion of the filtration generated by the Brownian motion. 
\par We consider the following backward stochastic differential equation (BSDE for short).

\begin{equation}\label{eq:11}
y_t=\xi+\int_t^T\!f(s,y_s,z_s)\,ds-\int_t^T\!z_s\,dW_s,\quad  t\in[0,T].
\end{equation}
where the generator $f:\Omega\times[0,T]\times\mathbb{R}\times\mathbb{R}^{1\times d}\rightarrow\mathbb{R}$ is a predictable function and terminal value $\xi$ is an $\mathcal{F}_T$-measurable random variable. 
\par The theory of BSDE is powerful to treat important issues arising in many applied fields such as finance and optimal control.  
A general nonlinear pricing problem of the European contingent claim in complete market is equivalent to solve the BSDE \eqref{eq:11}. 
In this case, $\xi$ is the contingent claim to hedge and $T$ is the maturity date.
Let us assume that \eqref{eq:11} has a solution $(y_t,z_t)$ in an appropriate space. If the generator is uniformly Lipschitz in $z$ (with Lipschitz constant $b$), we can apply the Girsanov measure change to the equation which leads to 

\begin{equation}\label{eq:12} 
y_t=\xi+\int_t^T\!f(s,y_s,0)\,ds-\int_t^T\!z_s\,dW_s^{\mathbb{Q}},\quad  t\in[0,T]. 
\end{equation}
where
 \[
\mathbb{Q}:=\exp\bigg(\int_0^T\!g(s,y_s,z_s)\,dW_s-\frac{1}{2}\int_0^T\!g^2 (s,y_s,z_s)\,ds\bigg)\cdot\mathbb{P},
\]
\[g(s,y_s,z_s )=\frac{f(s,y_s,z_s )-f(s,y_s,0)}{|z_s|^2}z_s\mathbb{1}_{|z_s|\neq 0},
\] 
and $W^\mathbb{Q}:=W-\int_0^{\cdot}\!g(s,y_s,z_s )\,ds$ is a $Q-$Brownian motion.

\par In finance, $\mathbb{Q}$ is called risk-neutral measure or martingale measure (see \cite{EPQ}). For convenience, let us assume that $f$ only depends on $z$ (hence $f=f(s,z)$ and $f(s,0)=0$. Then we have 
\[
y_t=\xi-\int_t^T\!z_s\,dW_s^{\mathbb{Q}}.
\]
When $\xi$ is square-integrable, it is well-known that the fair price of $\xi$ is evaluated as the expectation of the claim under $\mathbb{Q}$ (see e.g. \cite{EPQ}), that is,

\begin{equation}\label{eq:13}
y_t=\mathbb{E}^{\mathbb{Q}} [\xi |{\mathcal{F}_t}]. 
\end{equation}

At this point, one may be interested in looking for an "optimal" integrability condition, under which it is possible to represent the price by the risk-neutral measure, on terminal value. The paper of \cite{AIP} gave a partial resolution to this problem. Motivated by the expression \eqref{eq:13}, they introduced  the notation of measure solution which is benefit to give an efficient formula of pricing contingent claim by martingale measure. In Lipschitz setting, they showed the existence of the measure solution when the terminal value is $L^p$-integrable for $p>1$. In this case, one can use the H$\ddot{\text{o}}$lder's inequality and the boundness of moments of the exponential martingale to show $\mathbb{E}^{\mathbb{Q}}[\xi]<\infty$. If the terminal value is assumed to be integrable (i.e. $L^1$-integrable), it is not guaranteed that $\mathbb{E}^{\mathbb{Q}}[\xi]<\infty$, so the measure solution does not exist in general. That is, we need a stronger integrability condition on terminal value. Consequently, we want to find a sufficient integrability condition which is weaker than $L^p$-integrability for any $p>1$ and is stronger than $L^1$-integrability. 
\par Obviously, the expression \eqref{eq:13} is significant if and only if the following condition holds.

\begin{equation}\label{eq:14}
\mathbb{E}^{\mathbb{Q}}[|\xi|]=\mathbb{E}\biggl [|\xi|\exp (\int_0^T\!g(s,z_s )\,dW_s -\frac{1}{2}\int_0^T\!g^2 (s,z_s )\,ds ) \biggr ]<\infty.
\end{equation}

As $|g(s,z_s)|\le b$, above condition is equivalent to
\[
\mathbb{E}\biggl [|\xi|\exp\bigg(\int_0^T\!g(s,z_s)\,dW_s\bigg)\biggr ]<\infty.
\]

\cite{HT} showed the following useful inequalities.
\begin{itemize}
\item $e^x y\le e^{\frac{x^2}{2\mu^2}}+e^{2\mu^2}y\exp\big(\mu\sqrt{2\log(1+y)}\big)$.
\item $\mathbb{E}\biggl[\exp \big(\frac{1}{2\mu^2}|\int_0^T\!q_s\,dW_s|^2\big)\biggr ]\le [1-\frac{b^2}{\mu^2}T]^{-1/2}$ if $|q_s|\le b, \mu>b\sqrt{T}$.
\end{itemize}
From these two inequalities, we can deduce 
\[
\mathbb{E}\biggl[|\xi|\exp\bigg(\int_0^T \!g(s,z_s )\,dW_s\bigg)\biggr]\le\bigg [1-\frac{b^2}{\mu^2}T\bigg]^{-1/2}+e^{2\mu^2}\mathbb{E}\big[|\xi|\exp\big(\mu\sqrt{2\log(1+|\xi|)}\big)\big].
\]
So, we can get one sufficient condition to guarantee \eqref{eq:14} such that 
\[
\mathbb{E}\big[|\xi|\exp\big(\mu\sqrt{2\log(1+|\xi|)}\big)\big]<\infty.
\]
That is, $\xi$ is required to be $L\exp(\mu\sqrt{2\log(1+L)})$-integrable. Furthermore, if the condition \eqref{eq:14} is true, then $\xi$ is integrable under the measure $\mathbb{Q}$, so the BSDE \eqref{eq:11} is transferred into the BSDE \eqref{eq:12} with integrable parameters whose solution is called the $L^1$-solution. Also, the generating function of the equation \eqref{eq:12} does not depend on $z$, so the additional assumption (see \eqref{eq:15}) which is needed in the study of $L^1-$solution can be eliminated.  
\par \cite{PP} first introduced the notion of nonlinear BSDE and studied $L^2$-solution under the Lipschitz condition on generator.
\par \cite{BDHPS} studied  $L^p$-solutions ($p\ge1$) of BSDEs with monotonic generators. On the other hand, they introduced the following sub-linear growth assumption on generator to ensure the wellposedness of $L^1$-solution (hence $p=1$). 
\begin{equation}\label{eq:15}
|f(t,y,z)-f(t,0,0)|\le a|y|+b|z|^q,\quad (t,y,z)\in [0,T]\times\mathbb{R}\times\mathbb{R}^{1\times d} 
\end{equation}
for some $q\in [0,1)$. Later, \cite{F1} studied the wellposedness and comparisons of $L^p$-solutions ($p>1$) under various kinds of extended monotonicity conditions.
Also, \cite{F2} showed the existence, uniqueness and stability of $L^1$-solutions to BSDEs under one-sided Osgood condition, one of extended monotonicity conditions. However, one cannot find any results about the comparison principle of $L^1-$solutions. 
\par Recently, \cite{HT} studied the solution to BSDE in $L\exp(\mu\sqrt{2\log(1+L)})$-setting such that $\mu>\mu_0$ for some critical value $\mu_0$, that is, the terminal value is assumed to be $L\exp(\mu\sqrt{2\log(1+L)})$-integrable. This integrability is stronger than $L\log L$-integrability and weaker than $L^p$-integrability for any $p>1$. They showed the existence of solution to that BSDE under the linear growth condition on the generator. Furthermore they gave counterpart examples which show that $L\log L$-integrability is not sufficient to guarantee the existence of the solution. 
\par Afterwards, \cite{BHT} improved the existence result and gave the uniqueness result for the preceding BSDE under the Lipschitz condition by investigating the nice property of the solution $Y$ that $\phi(|Y|,\mu)$ belongs to class (D) (this nice property will be used effectively in our discussion). 
In their proof of uniqueness, the Lipschitz assumption played a crucial role because the representation of a solution to the linear BSDE was used. \cite{F3} studied the critical case: $\mu=\mu_0$. Note that if $\mu<\mu_0$, then the BSDE does not admit a solution in general (see \cite{HT}). In this paper, we state the uniqueness result under One-Sided Osgood condition, the extended form of monotonicity conditions. The next subject of this paper is to state the comparison principles. As it is well known, the comparisons for BSDEs are fundamental in the theory of nonlinear expectations, particularly in constructing the dynamic risk measures.  
\cite{CEP} showed a general comparison theorem by means of the super-martingale measure which is corresponded to the "no-arbitrage" condition in financial sense. In their paper, the terminal value was only assumed to guarantee the existence of a solution and the existence of certain super-martingale measure was also assumed, independently. 
For the BSDEs with $L\exp(\mu\sqrt{2\log(1+L)})$-integrable terminal values, we show the existence of such super-martingale measure. Then we can use directly the comparison theorems established by \cite{CEP}. This will just provide various applications to the world of dynamic risk measures in the same way as in \cite{CEP}. Also, we show the comparison theorem for BSDEs under one-sided Osgood condition (not Lipschitz in $y$) using penalization method. As the last subject, we state the stability result for the BSDEs with generators which is linear with respect to $z$ under One-Sided Osgood condition. 
The basic idea in all the proof in this paper is to associate the solution of the main BSDE with the $L^1$-solution of a certain BSDE with integrable parameters using Girsanov change, effectively.

\section{Notations and Assumptions}
\begin{itemize}
\item For $A\in\mathcal{F}$, $\mathcal{F}$-measurable random variable $\eta$ and probability measure $\mathbb Q$, we define $\mathbb{E}^{\mathbb{Q}}[\eta;A]:=\int_A\!\eta\,d{\mathbb{Q}}$. And $\mathbb{E}^{\mathbb{Q}}[\eta]:=\mathbb{E}^{\mathbb{Q}}[\eta;\mathbb{Q}]$. 

\item $\mathbb{T}(0,T)$ is a set of stopping times $\tau$ such that $0\le\tau\le T$.

\item For any predictable process $\phi$, we put $\mathcal{E}(\phi\bullet W):=\exp(\int_0^{\cdot}\!\phi_r\,dW_r-\frac{1}{2}\int_0^{\cdot}\!\phi_r^2\,dr)$.

\item We say that the process $Y=\{Y_t\}_{0\le t\le T}$ belongs to class $(D)$ if the family $\{Y_{\tau},\tau\in\mathbb{T}(0,T)\}$ is uniformly integrable.

\item $|\cdot|$ means the standard Euclidean norm.

\item $M^p ([0,T],\mathbb{R}^{1\times d};\mathbb{Q})$ is the space of predictable processes $Z$ with values in $\mathbb{R}^{1\times d}$ such that 
\[
|Z|_{M^p}:=\mathbb{E}^{\mathbb{Q}}\biggl [(\int_0^T\!|Z_s |^2\,ds)^{p/2}\biggr ]^{1\wedge1/p}<\infty.
\] 
If $\mathbb{Q}=\mathbb{P}$, then we denote it by $M^p([0,T];\mathbb{R}^{1\times d})$.

\item $H_T^1(\mathbb{Q})$ is the space of real c\`adl\`ag, adapted processes $Y$ such that $\mathbb{E}^{\mathbb{Q}}
[\sup_{t\in[0,T]}|Y_t|]<\infty$. If $\mathbb{Q}:=\mathbb{P}$, then we use $H_T^1$.

\item The solution of \eqref{eq:11} is denoted by a pair $\{(Y_t,Z_t ),t\in[0,T]\}$ of predictable processes with values in $\mathbb{R}\times\mathbb{R}^{1\times d}$ such that $Y$ is $\mathbb{P}$-a.s. continuous, $Z\in M^2 ([0,T];\mathbb{R}^{1\times d})$ and $(Y,Z)$ satisfies the equation \eqref{eq:11}.

\item For any real valued function $g$, we define $g^+ := \max (g, 0)$.
\end{itemize}

Define the real function $\psi$:
\[
\psi(x,\mu):=x\exp(\mu\sqrt{2\log(1+x)}),\quad (x,\mu)\in [0,\infty)\times(0,+\infty).
\]
Then, it has the following properties (see \cite{BHT, HT}).
\begin{itemize}
\item For any $x\in\mathbb{R}$ and $y\ge 0$, we have
\begin{equation}\label{eq:21}
e^x y\le e^{\frac{x^2}{2\mu^2}}+e^{2\mu^2}\psi(y,\mu).
\end{equation}	

\item Let $\mu>b\sqrt{T}$. Then for any $d$-dimensional adapted process $q$ with $|q_t|\le b$ a.s., for any $t\in[0,T]$, 
\begin{equation}\label{eq:22}
\mathbb{E}\biggl [\exp\bigg(\frac{1}{2\mu^2}\ \big|\int_t^T\!q_s\,dW_s\big|^2\bigg)\bigg{|}\mathcal{F}_t\biggr ]\le \big[1-\frac{b^2}{\mu^2}(T-t)\big]^{-1/2}.
\end{equation}

\item For any $\mu>0$, $\psi(\cdot,\mu)$ is convex, that is, for any $0\le \lambda\le1$ and $x,y\in[0,+\infty)$,
\begin{equation}\label{eq:23}
\psi(\lambda x+(1-\lambda)y,\mu)\le\lambda\psi(x,\mu)+(1-\lambda)\psi(y,\mu).
\end{equation}	

\item For any $l>1,x\le0$, we have
\begin{equation}\label{eq:24}
\psi(lx,\mu)\le\psi(l,\mu)\psi(x,\mu).	
\end{equation}
\end{itemize}

We present some useful assumptions on generator below.

\vskip 0.2cm
\textbf{(A1)} $f$ satisfies the One-Sided Osgood condition with respect to $y$, that is, there exists a non-decreasing and concave function $\rho:\mathbb{R}^+\rightarrow \mathbb{R}^+$ with $\rho(0)=0, \rho(t)>0$ for $t>0$ and $\int_{\mathbb{R}^+}\!\frac{dt}{\rho(t)}=+\infty$ such that for any $y,y'\in\mathbb{R}$ and $z\in\mathbb{R}^{1\times d}$,
\[
\frac{y-y'}{|y-y'|}\mathbb{1}_{|y-y'|\neq 0} (f(t,y,z)-f(t,y',z))\le \rho(|y-y'|). 
\]
\vskip 0.2cm
\textbf{(A2)} $f$ is uniformly Lipschitz in $z$, that is, there exists a constant $b$ such that for any $y\in\mathbb{R}$ and $z,z'\in\mathbb{R}^{1\times d}$,
\[
|f(t,y,z)-f(t,y,z')|\le b|z-z'|.
\]
\vskip 0.2cm
\textbf{(A3)} The map $y\mapsto f(t,y,z)$ is continuous.
\vskip 0.2cm
\textbf{(A4)} $f$ has linear growth in $y$, that is, there exists a constant $a\ge 0$ such that for 
any $y,y'\in\mathbb{R}$ and $z\in\mathbb{R}^{1\times d}$,
\[
|f(t,y,z)-f(t,0,z)|\le a|y|
\]
\vskip 0.2cm
\textbf{(A5)} $f$ is uniformly Lipschitz in $y$, that is, there exists a constant $r$ such that for any $y,y'\in\mathbb{R}$ and $z\in\mathbb{R}^{1\times d}$,
\[
|f(t,y,z)-f(t,y',z)|\le r|y-y'|.
\]

%%%%%%%%%%%%%%%%%%%%%%%%%%%%%%%%%%%%%
\section{Uniqueness}
\begin{theorem}\label{th:31}
Let assumptions \textbf{(A1), (A2)} hold. 
Then, BSDE \eqref{eq:11} has at most one solution $(Y,Z)$ such that $\psi(Y,c)$ belongs to class $(D)$ for some $c>0$.
\end{theorem}

\begin{proof} For $i=1,2$, let $(Y^i,Z^i)$ be a solution to \eqref{eq:11} such that $\psi(Y^i,c^i)$ belongs to the class $(D)$ for some $c^i>0$. Since $\psi(x,\mu)$ is non-decreasing in $\mu$, both $\psi(Y^1,c)$ and $\psi(Y^2,c)$ belong to class $(D)$ for $c=c^1\wedge c^2$. 

Define $(\bar Y,\bar Z):=(Y^1-Y^2,Z^1-Z^2)$.
For any $\tau\in\mathbb{T}(0,T)$, by \eqref{eq:23} and \eqref{eq:24},

\begin{align}
\psi(|\bar Y_{\tau}|,c)&\le \psi(|Y_{\tau}^1|+|Y_{\tau}^2|,c)=\frac{1}{2} \psi(2|Y_{\tau}^1|,c)+\frac{1}{2}\psi(2|Y_{\tau}^2 |,c)\nonumber\\
&\le\frac{1}{2} \psi(2,c)[\psi(|Y_{\tau}^1|,c)+\psi(|Y_{\tau}^2|,c)]\nonumber.
\end{align}

So, $\psi(|\bar Y|,c)$ is also belongs to class $(D)$.
\par We first restrict our discussion to the case of $T<\frac{c^2}{b^2}$ (hence $c>b\sqrt{T}$).  
Obviously, $(\bar Y,\bar Z)$ satisfies the following equation.
\begin{equation}\label{eq:31}
\bar Y_t=\int_t^T\!\bar f(s,\bar Y_s,\bar Z_s )\,ds-\int_t^T\! \bar Z_s\,dW_s,\quad  t\in[0,T].
\end{equation}
where $\bar f(s,y,z):=f(s,y+Y_s^2,z+Z_s^2 )-f(s,Y_s^2,Z_s^2)$.

We define
\[
\bar g(s,y,z):=\mathbb{1}_{|z|\neq 0}\frac{\bar f(s,y,z)-\bar f(s,y,0)}{|z|^2}z.  
\]
Then, it holds that
\[
\bar g_s:=\bar g(s,\bar Y_s,\bar Z_s )=\mathbb{1}_{|\bar Z_s|\neq 0}\frac{f(s,Y_s^1,Z_s^1)-f(s,Y_s^1,Z_s^2)}{|\bar Z_s|^2} \bar Z_s.
\]
From the assumption \textbf{(A2)}, we get $|\bar g|\le b$, a.s. and so $\mathbb{E}[\mathcal{E}(\bar g\bullet W)_t]=
\exp\big(\int_0^t\!\bar g_s \,dW_s -\frac{1}{2} \int_0^t\!\bar g_s^2\,ds\big)$ is an uniformly integrable martingale.  
\par By the virtue of Girsanov change, we have
\[
 \bar Y_t=\int_t^T\!\bar f(s,\bar Y_s,0)\,ds-\int_t^T\!\bar Z_s\,dW_s^{\mathbb{Q}},\quad  t\in[0,T].		
\]
where $\mathbb{Q}:=\mathcal{E}(\bar g_s\bullet W)_T\cdot\mathbb{P}$, and $W^{\mathbb{Q}}:=W-\int_0^{\cdot}\!\bar g_s\,ds$. 

\par Note that $\mathbb{Q}$ is a probability measure equivalent to $\mathbb{P}$ and $W^{\mathbb{Q}}$ is a Brownian motion under $\mathbb{Q}$. 
Then, for any $\tau\in\mathbb{T}(0,T)$ and $A\in\mathcal{F}$, by \eqref{eq:21} and \eqref{eq:22},
\begin{align}
\mathbb{E}^{\mathbb{Q}} (|\bar Y_{\tau}|;A)&=\mathbb{E}(|\bar Y_{\tau}|\cdot\mathcal{E}(\bar g\bullet W)_{\tau};A)\nonumber\\
&\le \mathbb{E}\biggl [|\bar Y_{\tau}|\exp\bigg(\int_0^T\!\bar g_s\,dW_s\bigg);A\biggr ]\le
\mathbb{E}\biggl [\exp\bigg(\frac{|\int_0^T\!\bar g_s\,dW_s|^2}{2c^2 }\bigg)+e^{2c^2} \psi(|\bar Y_{\tau}|,c);A\biggr ]\nonumber\\
&\le (1-\frac{b^2}{c^2}T)^{-1/2}+\mathbb{E}[e^{2c^2} \psi(|\bar Y_{\tau}|,c);A]\nonumber.
\end{align}
So, $\bar Y$ belongs to class $(D)$ under $\mathbb{Q}$. 
Now we give an estimate on $\bar Z$ under $\mathbb{Q}$.
\par Let $1<p<\infty$, $p^{-1}+q^{-1}=1$ and $k=\frac{\sqrt{p}}{2(\sqrt{p}-1)}$. 
Then for any $\tau\in\mathbb{T}(0,T)$,
\[
\mathcal{E}((k\bar g)\bullet W)_{\tau}=\exp\big(\int_0^\tau\!k\bar g_s \,dW_s \big)\exp\big(-\frac{1}{2} \int_0^\tau\!k^2 \bar g_s^2\,ds\big)\ge\exp\big(\int_0^{\tau}\!k\bar g_s\,dW_s\big)\exp\big(-\frac{1}{2} b^2 k^2 T\big).
\]
Since $\mathbb{E}[\mathcal{E}((k\bar g )\bullet W)_{\tau}]=1$, we obtain
\[
\sup_{\tau\in\mathbb{T}(0,T)}\mathbb{E}[\exp(\int_0^{\tau}\!k\bar g_s\,dW_s)]\le\exp(\frac{1}{2} b^2 k^2 T).
\]
Therefore, according to \cite{K}, Theorem 1.5, $\mathcal{E}(\bar g\bullet W)$ is $L^q$-bounded martingale. Using H$\ddot{\text{o}}$lder's inequality, we obtain
\begin{align}
\mathbb{E}^{\mathbb{Q}} \biggl [\bigg(\int_0^T\!|\bar Z_s|^2\,ds\bigg)^{1/p}\biggr ]&=\mathbb{E}\biggl [\bigg(\int_0^T\!|\bar Z_s |^2\,ds\bigg)^{1/p} \mathcal{E}(\bar g\bullet W)_T\biggr ]\nonumber\\
&\le\mathbb{E}\biggl [\int_0^T\!| \bar Z_s|^2 \,ds\biggr ]^{1/p}\cdot \mathbb{E} [\mathcal{E}(\bar g\bullet W)_T^q ]^{1/q}<\infty\nonumber. 
\end{align}

Taking $\bar p:=\frac{2}{p}$, then $\bar Z\in M^{\bar p}([0,T],\mathbb{R}^{1\times d};\mathbb{Q})$. Moreover, due to the arbitrariness of $p$, it holds that $\bar Z\in M^{\bar p} ([0,T],\mathbb{R}^{1\times d};\mathbb{Q})$ for any $0<\bar p<2$.

\par Therefore, $(\bar Y,\bar Z)$ is an $(L^1-)$ solution of the following BSDE such that  $\bar Y$ belongs to class $(D)$ and $\bar Z\in M^{\bar p} (\mathbb{R}^{1\times d};\mathbb{Q})$ for any $0<\bar p<2$.
\begin{equation}\label{eq:32}
y_t=\int_t^T\!\bar f(s,y,0)\,ds-\int_t^T\!z_s\,dW_s^{\mathbb{Q}},\quad  t\in[0,T].
\end{equation}

Since $\bar f(s,0,0)=0$, a pair $(0,0)$ is also a solution of \eqref{eq:32}. 
\par On the other hand, for any $y,y'\in\mathbb{R}$,
\begin{align}
\frac{y-y'}{|y-y'|} \mathbb{1}_{|y-y'|\neq 0}& (\bar f(s,y,0)-\bar f(s,y',0))\nonumber\\
&=\frac{y-y'}{|y-y'|} \mathbb{1}_{|y-y'|\neq 0} (f(s,y+Y_s^2,Z_s^2 )-f(s,y'+Y_s^2,Z_s^2 ))\le \rho(|y-y'|)\nonumber.
\end{align}

Therefore, according to the uniqueness of $L^1$-solution of BSDEs with generators of One-Sided Osgood type (see \cite{F2}, Theorem 1), we have $(\bar Y,\bar Z)=(0,0)$. For larger value of $T$, we first discuss on interval $[T-\delta,T]$ for small $\delta>0$ from which we get $(\bar Y_t,\bar Z_t)=0$ for $T-\delta\le t\le T$ and then with the terminal value $\bar Y_{T-\delta}=0$, we discuss on interval $[T-2\delta,T-\delta]$ from which $(\bar Y_t,\bar Z_t)=0$ for $T-2\delta\le t\le T-\delta$ and so on by an inductive argument. This provides $(\bar Y,\bar Z)=0$ on the whole interval $[0,T]$. That is, we have $Y^1=Y^2$ and $Z^1=Z^2$.
\end{proof}

Due to the existence result (\cite{BHT}, Theorem 2.4), we get the following result.
\begin{corollary}\label{co:32}
Suppose that \textbf{(A1), (A2), (A3)} and \textbf{(A4)} hold. We further assume that there exists a constant  $\mu>b\sqrt{T}$ such that 

\begin{equation}\label{eq:33}
\psi\bigg(|\xi|+\int_0^T\!|f(t,0,0)|\,dt,\mu\bigg)\in L^1(\Omega,\mathbb{P}).
\end{equation}

Then, BSDE \eqref{eq:11} has a unique solution $(Y,Z)$ such that $\psi(Y,c)$ belongs to the class $(D)$ for some $c>0$. Moreover we have the following estimate on $Y$.
\begin{equation}\label{eq:34}
|Y_t |\le \frac{1}{\sqrt{1-\frac{b^2}{\mu^2}(T-t)}}e^{a(T-t)} +e^{2\mu^2+a(T-t)} \mathbb{E}\biggl [\psi\bigg(|\xi|+\int_t^T\!|f(s,0,0)|\,ds,\mu\bigg)\bigg{|}\mathcal{F}_t\biggr ].
\end{equation}
\end{corollary}
	
\begin{remark}\label{re:33}
We also have an estimate on $\psi(Y,\mu)$ (See the last inequality in the proof of \cite{BHT}, Theorem 2.4). For some constants $A,B\ge 0$, it holds that
\begin{equation}\label{eq:35}
\psi(|Y_t |,\mu)\le A+B\cdot\mathbb{E} \biggl [\psi\bigg(|\xi|+\int_0^T\!|f(s,0,0)|\,ds,\mu\bigg)\bigg{|}\mathcal{F}_t\biggr ].
\end{equation}
\end{remark}
%%%%%%%%%%%%%%%%%%%%%%%%%%%%%%%%%%%%%%%%%%%%%%%%%%%%%%%%%%%
\section{Comparisons}
We first show the comparison principle for the BSDE with Lipschitz generator.

\begin{theorem}\label{th:41}
Let $(\xi,f)$ and $(\xi',f')$ be any two pairs of terminal value and generator of \eqref{eq:11}, respectively. Let $(Y,Z)$ and $(Y',Z')$ be associated solutions such that $\psi(Y,c)$ and $\psi(Y',c')$ belong to class $(D)$ for some $c,c'>0$. Suppose that $f$ satisfies $\textbf{(A2)}$ and $\textbf{(A5)}$.  If $\xi\le\xi'$ and $f(t,Y',Z')\le f'(t,Y',Z')$ then $Y_t\le Y'_t$ for all $t\in[0,T]$, $\mathbb{P}$-a.s. Moreover this comparison is strict, that is, if $Y_t^1=Y_t^2,\mathbb{P}-a.s.$ on $A\in\mathcal{F}_t$, then $Y_s^1=Y_s^2$ on $[t,T]\times A$ up to evanescence.
\end{theorem}

\begin{proof} As we showed at the beginning part of the proof of theorem \ref{th:31}, $\psi(Y-Y',c_0)$ belongs to class $(D)$ for $c_0=c\wedge c'$.
We assume that $c_0>b\sqrt{T}$ without loss of generality. For larger $T$, we can adopt the same strategy as in the proof of theorem \ref{th:31}.
~\\
Let us define the process:
\[
\Gamma_s:=\frac{f(s,Y'_s,Z_s )-f(s,Y'_s,Z'_s )}{|Z_s-Z'_s |^2}\mathbb{1}_{|Z_s-Z'_s |=0}(Z_s-Z'_s)
\]
 which is uniformly bounded. The measure $\mathbb{Q}$ is defined as follows.
\[ 
\frac{d\mathbb{Q}}{d\mathbb{P}}:=\mathcal{E}(\Gamma\bullet W)_T=\int_0^T\!\Gamma_s\,dW_s -\frac{1}{2} \int_0^T\!\Gamma_s^2 \,ds.
\]
Then, $W^{\mathbb{Q}}:=W-\int_0^{\cdot}\!\Gamma_s\,ds$ is $\mathbb{Q}$-Brownian motion. As we showed in preceding discussion, $Z-Z'\in M^q (\mathbb R^{1\times d};\mathbb{Q})$ for any $q\in(1,2)$.
Therefore, 
\[
-\int_0^t\!(f(s,Y'_s,Z_s )-f(s,Y'_s,Z'_s ))\,ds+\int_0^t\!(Z_s-Z'_s)d\,W_s =\int_0^t\!(Z_s-Z'_s)\,dW_s^{\mathbb{Q}}.
\]
is a $\mathbb{Q}$-martingale. On the other hand, 
\begin{align}
\mathbb{E}^{\mathbb{Q}}\biggl [\int_0^T\!(f(s,Y_s,Z_s )-&f(s,Y'_s,Z_s ))\,ds\biggr ]\le r\mathbb{E}^{\mathbb{Q}} \biggl [\int_0^T\!|Y_s-Y'_s|\,ds\biggr ]\nonumber\\
&=r\mathbb{E}\biggl [\int_0^T\!|Y_s-Y'_s|\mathcal{E}(\Gamma\bullet W)_T\,ds\biggr ]\nonumber\\
&\le r\mathbb{E}\biggl [\int_0^T\!\bigg(1-\frac{b^2}{c_0^2} T\bigg)^{-1/2}+e^{2c_0^2} \psi(|Y_s-Y'_s |,c_0 )\,ds\biggr ]\nonumber\\
&\le rT\bigg(1-\frac{b^2}{c_0^2}T\bigg)^{-1/2}+rTe^{2c_0^2 }\cdot 
\sup_{\tau\in\mathbb{T}}\mathbb{E}[\psi(|Y_{\tau}-Y'_{\tau}|,c_0)]<\infty\nonumber.
\end{align}
Now, both comparison and strict comparison theorems just follow from \cite{CEP}, Theorems 1, 2 and 3.
\end{proof}

\begin{remark}\label{re:42}
As an immediate consequence of Theorem \ref{th:41}, we can see that the solution of the BSDE \eqref{eq:11} is unique. So, we have provided an alternative method for the proof of the uniqueness part in Lipschitz setting than that of \cite{BHT}.
\end{remark}

Now we discuss the comparison theorem under one-sided Osgood condition.

\begin{theorem}\label{th:43}
The comparison theorem still holds under the assumptions \textbf{(A1)} and \textbf{(A2)}.
\end{theorem}

\begin{proof} Set $\bar{\xi}:=\xi-\xi', \bar f:=f-f', \delta f(s,y,z):=f(s,y+Y'_s,z+Z'_s)-f(s,Y'_s,Z'_s)$.  
The pair $(\bar Y,\bar Z ):=(Y-Y',Z-Z')$ satisfies
\[
\bar Y_t= \bar{\xi}+\int_t^T\![\delta f(s,\bar Y_s,\bar Z_s )+\bar f(s,Y'_s,Z'_s)]\,ds-\int_t^T\!\bar Z_s\,dW_s.
\]
After an application of the Girsanov change, we have
\begin{equation}\label{eq:41}
\bar Y_t= \bar{\xi}+\int_t^T\![\delta f(s,\bar Y_s,0)+\bar f(s,Y'_s,Z'_s)]\,ds-\int_t^T\!\bar Z_s\,dW_s^{\mathbb{Q}},
\end{equation}
where the probability measure $\mathbb{Q}$ is similarly defined as before.
\par Note that $\mathbb{E}^{\mathbb{Q}} [\bar{\xi}]<\infty,\ \delta f(s,0,0)=0$ and $\bar Y_t$ belongs to class $(D)$ under $\mathbb{Q}$. We also note that $\bar Z\in M^q (\mathbb R^{1\times d};\mathbb Q)$ for any $q\in (1,2)$.
Applying Tanaka's formula to \eqref{eq:41},
\begin{equation}\label{eq:42}
\bar Y_{t}^+=\bar{\xi}^+ +\int_{t}^{T}\!\mathbb{1}_{\bar Y_s>0}[\delta f(s,\bar Y_s,0)+\bar f(s,Y'_s,Z'_s)]\,ds -\int_{t}^{T}\!\mathbb{1} _{\bar Y_s >0} \bar Z_s\,dW_s^{\mathbb{Q}} -\frac{1}{2} L_t^0,
\end{equation}
where $L_t^0$ is the local time of $\bar Y_t$ at $0$, it is an increasing process such that $L_0^0=0$.
Since $\bar f(s,Y'_s,Z'_s )\le 0$, we see that
\[
\mathbb{1}_{\bar Y_s>0}[\delta f(s, \bar Y_s,0)+\bar f(s,Y'_s,Z'_s )]\le \mathbb{1}_{\bar Y_s>0}\cdot \delta f(s, \bar Y_s,0)=\mathbb{1}_{\bar Y_s>0}\frac{\bar Y_s}{|\bar Y_s|}\cdot \delta f(s, \bar Y_s,0)\le\rho(\bar Y_s^+).
\]
On the other hand, the function $\rho(\cdot)$ has linear growth since it is non-decreasing and concave valued $0$ at $0$. If we denote by $l$ the linear growth, then $\mathbb E^{\mathbb Q}[\rho(\bar Y_s^+)]\le \mathbb E^{\mathbb Q}[l(\bar Y_s^++1)]=l\mathbb (\mathbb E^{\mathbb Q}[\bar Y_s^+]+1)<\infty$. 
~\\
Taking conditional expectations on both sides of \eqref{eq:42} with respect to $\mathbb Q$, we get
\[
\mathbb E^{\mathbb Q}[\bar Y_t^+|\mathcal F_t] \le \int_t^T\!\mathbb{E}^{\mathbb{Q}} [\rho(\bar Y_s^+)|\mathcal F_t]\,ds\le \int_t^T\! \rho\big (\mathbb{E}^{\mathbb{Q}}[\bar Y_s^+|\mathcal F_t]\big )\,ds
\]
where we used Jensen's inequality and  $\bar\xi ^+=0$. Then, Bihari's inequality implies that $\bar Y_t^+=0,\ \mathbb Q-a.s.$ for each $t\in [0,T]$.
As $\mathbb{Q}$ is equivalent to $\mathbb{P}$, we have $\bar Y_t^+=0,\ \mathbb{P}-a.s.$ Hence $Y_t^1\le Y_t^2,\ \mathbb {P}-a.s.$
\end{proof}
\begin{remark}
Theorem \ref{th:43} can be regarded as a generalization of Theorem \ref{th:31}, as one sees easily.
In general, the strict comparison theorem does not hold in a monotonicity setting (see \cite{PR}, pp. 416).  
\end{remark}
%%%%%%%%%%%%%%%%%%%%%%%%%%%%%%%%%%%%%%%%%%%%%%%%%%%%%%%%%%%%%%%%
\section{Stability}
In this section, we state the stability result for BSDE \eqref{eq:11}. We shall restrict to the case where the generator is linear with respect to $z$. The more general case is left for the future work. 
Before we study the stability, we give the following useful lemma.
\begin{lemma}\label{le:51}
Suppose that the generator satisfies \textbf{(A1), (A2), (A3)} and \textbf{(A4)}. Instead of \eqref{eq:33}, we assume that 
\[
\sup_{t\in[0,T]}\mathbb{E}\biggl [\psi\bigg(|\xi|+\int_0^T\!|f(t,0,0)|dt,\mu\bigg)\bigg{|} \mathcal{F}_t\biggr ]\in L^1 (\Omega,\mathbb{P}),\mu>b\sqrt{T}.
\]
Then the BSDE \eqref{eq:11} has a unique solution such that $\psi(|Y|,\mu)\in H_T^1$.
\end{lemma}
\begin{proof} By Corollary \ref{co:32}, \eqref{eq:11} has a unique solution $(Y,Z)$ such that $\psi(Y,\mu)$ belongs to class $(D)$. Due to \eqref{eq:35}, we can see that $\psi(|Y|,\mu)\in H_T^1$.
\end{proof}

\begin{theorem}\label{th:52}
For each $n\in\mathbb{N}_0$, let us consider the following BSDEs depending on parameter $n$:
\[
Y_t^n=\xi^n+\int_t^T\!f^n (s,Y_s^n,Z_s^n)\,ds-\int_t^T\!Z_s^n\,dW_s,\quad t\in [0,T].
\]
We introduce the following assumptions.

\begin{enumerate}
\item For all $n$, $\xi^n$ and $f^n$ satisfy \textbf{(A1), (A2), (A3)} and \textbf{(A4)} with the same parameters $\rho(\cdot),a,b$.
\item $f^0$ is linear with respect to $z$, that is, $f^0(s,y,z)=f^0(s,y,0)+bz$. 
\item There exists a constant $\mu>b\sqrt{T}$ such that
\[
\psi\bigg(\xi^0+\int_0^T\!f^0 (t,0,0)\,dt,\mu\bigg)\in L^1(\Omega,\mathbb{P}).
\]
\item There exists a non-negative real sequence $(l_n )_{n=1,2,...}$ which converges to $0$ such that for each $n$, for any $(y,z)\in \mathbb{R}\times \mathbb{R}^{1\times d}$,
\[
|f^n (s,y,z)-f^0 (s,y,z)|\le l_n,\quad d\mathbb{P}\times dt-a.s.
\]
\item There exists a random variable $\eta$ satisfying $\psi(\eta,\mu)\in L^1(\Omega,\mathbb{P})$ such that
$|\xi^n-\xi^0|\le\eta$ for any $n\ge 1$ and $\mathbb{E}[|\xi^n-\xi^0|]\rightarrow 0, \quad n\rightarrow \infty$. 
\item There exists a constant $\mu>b\sqrt{T}$ such that
\[
\sup_{t\in[0,T]}\mathbb{E}\biggl [\psi\bigg(|\xi^0|+\int_0^T\!|f^0(t,0,0)\,dt,\mu\bigg)\bigg{|} \mathcal{F}_t\biggr ]\in L^1 (\Omega,\mathbb{P}).
\]
\item There exists a random variable $\eta$ satisfying $\psi(\sup_{t\in [0,T]}\mathbb{E}[\eta |\mathcal{F}_t],\mu)\in L^1(\Omega,\mathbb{P})$ such that $|\xi^n-\xi^0|\le\eta$ for any $n\ge 1$ and
\[
\mathbb{E}\bigg[\sup_{t\in [0,T]}\mathbb{E}\big[|\xi^n-\xi^0|\big|\mathcal{F}_t\big]\bigg]\rightarrow 0,\quad n\rightarrow\infty.
\]
\end{enumerate}
\begin{enumerate}
\item[(i)] Under assumptions 1-5, we have  
\[
\sup_{t\in[0,T]}\mathbb{E}\big[\psi(|Y_t^n-Y_t^0|,\mu)\big]\rightarrow 0, n\rightarrow\infty.
\] 
and for any $\beta\in (0,1)$,
\[
\mathbb{E}\biggl [\sup_{t\in[0,T]}\psi(|Y_t^n-Y_t^0|,\mu)^{\beta}+\bigg(\int_0^T\!|Z_s^n-Z_s^0|^2\, ds\bigg)^{\beta/2}\biggr ]\rightarrow 0, n\rightarrow\infty.
\]

\item[(ii)] Moreover, if assumptions 3,5 are replaced by assumptions 6,7, then it holds that 
\[
\mathbb{E}\bigl [\sup_{t\in[0,T]}\psi(|Y_t^n-Y_t^0|,\mu)\bigr ] \rightarrow 0,\ n\rightarrow\infty.
\]
\end{enumerate}
\end{theorem}
\begin{remark}
If assumptions 3-5 (resp., assumptions 4,6,7) are true, then  it just follows from the expressions \eqref{eq:23} and \eqref{eq:24} that $\psi\big(|\xi^n|+\int_0^T\!|f^n(t,0,0)|\,dt,\mu\big)\in L^1 (\Omega,\mathbb{P}).$ (resp., $\sup_{t\in[0,T]}\mathbb{E}\big [\psi\big(|\xi^n|+\int_0^T\!|f^n(t,0,0)|\,dt,\mu\big)\big{|} \mathcal{F}_t\bigr ]\in L^1 (\Omega,\mathbb{P}).$) for each $n\in \mathbb N$. Note that assumptions 6,7 are stronger than assumptions 3,5, respectively.
\end{remark}
\begin{proof}[\textbf{Proof of Theorem 5.2}]
~\\
$\text{(i).}$
By the virtue of Girsanov change, we have for each $n\in \mathbb{N}_0$,
\[
Y_t^{n}=\xi^{n}+\int_t^T\!f^{n} (s,Y_s^{n},0)\,ds-\int_t^T\!Z_s^{n}\,dW_s^{\mathbb{Q}^n}.
\]
where
\[
\frac{d\mathbb{Q}^n}{d\mathbb{P}}:=\mathcal{E}(g^{n}\bullet W)_T, \quad W^{\mathbb{Q}^n}:=W-\int_0^{\cdot}\!g^{n} (s)\,ds. 
\]
\[
g^n(s):=\frac{f^n(s,Y_s^n,Z_s^n)-f^n(s,Y_s^n,0)}{|Z_s^n|^2}\mathbb{1}_{|Z_s^n|\neq 0}Z_s^n.
\]
We put $\mathbb{Q}:=\mathbb{Q}^0$. Clearly, $W^{\mathbb{Q}^n}=W^{\mathbb{Q}}-\int_0^{\cdot}\! (g^n(s)-g^0(s))\,ds$ for each $n$.
So, we get
\[
Y_t^{n}=\xi^{n}+\int_t^T \bar f^{n} (s,Y_s^{n},Z_s^{n})ds-\int_t^T\!Z_s^{n}\,dW_s^{\mathbb{Q}}.
\]
where $\bar f^n(s,y,z):=f^n(s,y,z)-g^0(s)z=f^n(s,y,z)-bz$.
~\\
Note that $\bar f^0(s,y,z)=f^0(s,y,0)$. 
The same arguments as in the proof of preceding results give that
\[
\mathbb{E}^{\mathbb{Q}}[\eta]<\infty, \quad \mathbb{E}^{\mathbb{Q}}\biggl [\xi^{n}+\int_0^T\!\bar f^{n}(s,0,0)\biggr ]=\mathbb{E}^{\mathbb{Q}}\biggl [\xi^{n}+\int_0^T\!f^{n}(s,0,0)\biggr ]<\infty. 
\]
Moreover, both processes $Y^n$ and $\psi(|Y^n|,\mu)$ belong to class (D) under $\mathbb{Q}$ and  $Z^n\in M^{\bar p} ([0,T],\mathbb{R}^{1\times d};\mathbb{Q})$, for any $0<\bar p<2$.
~\\ 
And $\bar f^n$ has the sublinear growth in $z$ from
\begin{align}
|\bar f^n(s,y,z)-\bar f^n(s,y,0)|&=|f^n(s,y,z)-f^n(s,y,0)-bz|\nonumber\\
&\le 2l_n+|f^0(s,y,z)-f^0(s,y,0)-bz|=2l_n.\nonumber
\end{align}
Therefore, for each $n$, $(Y^n,Z^n)$ is a unique $L^1-$solution of the following BSDE under $\mathbb{Q}$.
\begin{equation}\label{eq:51}
y_t=\xi^{n}+\int_t^T\! \bar f^{n} (s,y_s,z_s)\,ds-\int_t^T\!z_s\,dW_s^{\mathbb{Q}}.
\end{equation}
From the assumption, $\xi^n$ converges to $\xi^0$ in probability and so does under $\mathbb{Q}$. As $|\xi^n-\xi^0|\le\eta, n\in \mathbb{N}_0$ and $\mathbb{E}^{\mathbb{Q}}[\eta]<\infty$, by Lebesgue's dominated convergence theorem, we get $\mathbb{E}^{\mathbb{Q}}[|\xi^n-\xi^0|]=0$.
Also, it holds that $|\bar f^n (s,y,z)-\bar f^0 (s,y,z)|=|f^n(s,y,z)-f^0(s,y,z)|\le l_n$ for any $n\in\mathbb{N}_0$.
Now we can use the stability results of $L^1$-solutions to BSDE \eqref{eq:51}. According to \cite{F2}, Theorem 4, it holds that

\begin{equation}\label{eq:52}
\sup_{t\in [0,T]} \mathbb{E}^{\mathbb{Q}}[|Y_t^n-Y_t^0|] \rightarrow 0,\ n\rightarrow\infty.	
\end{equation}
and for any $\beta\in (0,1)$,

\begin{equation}\label{eq:53}
\mathbb{E}^{\mathbb{Q}} \biggl [\sup_{t\in[0,T]}|Y_t^n-Y_t^0|^{\beta} +\bigg(\int_0^T\!|Z_s^n-Z_s^0|^2\,ds\bigg)^{\beta/2} \biggr ]\rightarrow 0,\ n\rightarrow\infty.	
\end{equation}
~\\
For the simplicity, we define
\[
Y^{n,0}:=Y^n-Y^0, Z^{n,0}:=Z^n-Z^0,\xi^{n,0}:=\xi^n-\xi^0, \bar f^{n,0}:=\bar f^n-\bar f^0=f^n-f^0=:f^{n,0}.
\]
The expression \eqref{eq:52} implies that $Y_t^{n,0}\xrightarrow{\mathbb Q} 0$ uniformly in $t$. As the measure $\mathbb{Q}$ is equivalent to $\mathbb P$, we see that $Y_t^{n,0}\xrightarrow{\mathbb P} 0$ uniformly in $t$. Moreover, it follows that $\psi(|Y_t^{n,0}|,\mu)\xrightarrow{\mathbb P} 0$ uniformly in $t$ from the fact that $\psi(\cdot,\mu)$ is strictly increasing.
~\\
On the other hand, using the expression \eqref{eq:35},
\begin{align}\label{eq:54}
\psi(|Y_t^{n,0}|,\mu)&\le A+B\cdot\mathbb E\biggl[\psi\bigg(|\xi^{n,0}|+\int_0^T\!|f^{n,0} (s,0,0)|\,ds,\mu\bigg)\bigg{|}\mathcal F_t\biggr]\nonumber\\
&\le A+B\cdot\mathbb E\biggl[\psi\bigg(\eta+T\sup_{n}l_n,\mu\bigg)\bigg{|}\mathcal F_t\biggr]=:(*).
\end{align}
\[
\sup_{t\in[0,T]}\mathbb{E}[(*)]=\mathbb{E}[(*)]\le A+\frac{1}{2}B\psi(2,\mu)\cdot\big[\psi(\eta,\mu) +\psi(T\sup_{n}l_n,\mu)\big]<\infty.
\]
So, by Lebesgue's dominated convergence theorem, we obtain
\[
\sup_{t\in[0,T]}\mathbb{E}\big[\psi(|Y_t^{n,0}|,\mu)\big] \rightarrow0,\ n\rightarrow\infty.
\]
From the expression \eqref{eq:53}, the process $|Y_t^{n,0}|^{\beta}\xrightarrow{ucp}0$ under $\mathbb Q$, so does under $\mathbb P$. Since $\psi(\cdot,\mu)$ is strictly increasing, $\psi(|Y_t^{n,0}|^{\beta},\mu)\xrightarrow{ucp}0$.
Using \eqref{eq:54} and  \cite{BDHPS}, Lemma 6.1, we deduce for any $\beta\in(0,1)$,
\begin{align}
\mathbb E\big[\sup_n\sup_{t\in[0,T]}\psi(|Y_t^{n,0}|,\mu)^{\beta}\big]&\le \mathbb E\big[\sup_n\sup_{t\in[0,T]}(*)^{\beta}\big]=\mathbb E\big[\sup_{t\in[0,T]}(*)^{\beta}\big]\nonumber\\
&\le A^{\beta}+B^{\beta}\frac{1}{1-\beta}\mathbb E\big[\psi(\eta+T\sup_n l_n,\mu)\big]^{\beta}<\infty.\nonumber
\end{align}
Then Lebesgue's dominated convergence theorem ensures that 
\[
\mathbb{E}\big [\sup_{t\in[0,T]}\psi(|Y_t^{n,0}|,\mu)^{\beta}\big]\rightarrow 0, n\rightarrow\infty.
\]
Next, for any $\varepsilon<(1-\beta)/\beta$, by H$\ddot{\text{o}}$lder's inequality, 
\begin{align}
\mathbb{E}\biggl [\bigg(\int_0^T\!|Z_s^{n,0}|&^2\,ds\bigg)^{\beta/2}\biggr ]=\mathbb{E}\biggl [\bigg(\int_0^T\!|Z_s^{n,0}|^2\,ds\bigg)^{\beta/2} [\mathcal{E}(b\bullet W)_T ]^{1/(1+\varepsilon)}\cdot [\mathcal{E}(b\bullet W)_T]^{{-1}/(1+\varepsilon)}\biggr ]\nonumber\\
&\le \mathbb{E}\biggl [\bigg (\int_0^T\!|Z_s^{n,0}|^2 \,ds\bigg )^{\beta(1+\varepsilon)/2}\mathcal{E}(b\bullet W)_T \biggr]^{1/(1+\varepsilon)}\cdot\mathbb{E}\big([\mathcal{E}(b\bullet W)_T ]^{-1/\varepsilon}\big)^{\varepsilon/(1+\varepsilon)}\nonumber\\
&=\mathbb{E}^{\mathbb{Q}}\biggl [\bigg (\int_0^T\!|Z_s^{n,0}|^2\,ds\bigg )^{\beta(1+\varepsilon)/2}\biggr ]^{1/(1+\varepsilon)}\cdot\mathbb{E}\big([\mathcal{E}(b\bullet W)_T]^{-1/\varepsilon}\big)^{\varepsilon/(1+\varepsilon)}\nonumber.
\end{align}

By \eqref{eq:53}, the last term tends to $0$ as $n\rightarrow\infty$. Consequently, we have
\[
\mathbb{E}\biggl [\sup_{t\in[0,T]}\psi(|Y_t^{n,0}|,\mu)^{\beta} +\bigg (\int_0^T\!|Z_s^{n,0}|^2\,ds\bigg )^{\beta/2}\biggr ]\rightarrow 0,\ n\rightarrow\infty.
\]
which is the desired result.

$\text{(ii).}$ We can have very similar procedure as in the proof of the first assertion, so we only sketch the proof.
Due to assumption 7, we see that $\mathbb E^{\mathbb Q}\bigg[\sup_{t\in [0,T]}\mathbb{E}\big[|\xi^n-\xi^0|\big|\mathcal{F}_t\big]\bigg]\rightarrow 0$, as $n\rightarrow\infty$ by Lebesgue's dominated convergence theorem. 
Using lemma \ref{le:51}, we deduce that $Y\in H_T^1(\mathbb Q)$. Then, according to the stability result of $L^1-$solution (\cite{F2}, Theorem 5), it holds that
\[
\mathbb{E}^{\mathbb{Q}} \biggl [\sup_{t\in[0,T]}|Y_t^{n,0}| +\bigg(\int_0^T\!|Z_s^{n,0}|^2\,ds\bigg)^{1/2} \biggr ]\rightarrow 0,\ n\rightarrow\infty.
\]
It can be easily seen that $\psi(|Y_t^{n,0}|,\mu)\xrightarrow{ucp}0, n\rightarrow\infty$ under $\mathbb P$. 
By the expression \eqref{eq:54} and assumption 7,
\[
\mathbb E\big[\sup_n\sup_{t\in[0,T]}\psi(|Y_t^{n,0}|,\mu)\big]\le \mathbb E\big[\sup_n\sup_{t\in[0,T]}(*)\big]=\mathbb E\big[\sup_{t\in[0,T]}(*)\big]<\infty.
\]
Now, we can use dominated convergence theorem to get the conclusion.
\end{proof}
\begin{remark}
Perhaps, one can try to prove directly the stability theorem without using the properties of $L^1-$solution. But this is not the objective within our framework. 
\end{remark}
\begin{remark}
One can easily check that the framework of this paper is also adapted to the critical case of $\mu=b\sqrt{T}$ due to the counter existence result of solution (see \cite{F3}).  
\end{remark}
%%%%%%%%%%%%%%%%%%%%%%%%%%%%%%%%%%%%%%%%%%%%%%%%%%%%%%%%%%%%%%%%%%%
%%                                                               %%
%% Use the two commands below for producing your bibliography    %%
%% with bibtex, then comment again the commands and include the  %%
%% content of the .bbl file in this file below the commands.     %%
%%                                                               %%
%%%%%%%%%%%%%%%%%%%%%%%%%%%%%%%%%%%%%%%%%%%%%%%%%%%%%%%%%%%%%%%%%%%
%\bibliographystyle{amsplain}
%\bibliography{yourbibfilename}
% add below the content of your .bbl file produced by bibtex.

\end{document}